\documentclass[a4paper,11pt]{article}

\usepackage{amsfonts,qtree,stmaryrd,textcomp}
\usepackage{pifont,amssymb,amsmath,amsthm,tabularx,graphicx}
\usepackage{tree-dvips,pictexwd,dcpic}
\usepackage{bbm}

\usepackage{stmaryrd}

\usepackage[T1]{fontenc}
\usepackage[latin1]{inputenc}
\usepackage[text={17cm,26cm},centering]{geometry}

\usepackage{etex}
\usepackage{fancyhdr}
\usepackage{graphicx}
\usepackage{enumitem}
\usepackage{amsmath}
\usepackage[bbgreekl]{mathbbol}
\usepackage{amssymb}
\usepackage{amsthm}
\usepackage[all]{xy}

\usepackage{epigraph}
\RequirePackage{bussproofs}

\usepackage{framed}
\usepackage{mathtools}
\usepackage{url}
\usepackage{proof}
\usepackage{xspace}
\usepackage{listings}
\usepackage{wrapfig}
\usepackage{tikz}
\usepackage{tikz-cd}

\usepackage{relsize}

\usetikzlibrary{positioning}
\usetikzlibrary{shapes}

\usetikzlibrary{arrows}
\usetikzlibrary{calc}
\usepackage{tikz-3dplot}
\usetikzlibrary{decorations.pathmorphing}

\usepackage{amsfonts,amssymb,amsmath,qtree,stmaryrd,textcomp, amsthm}

\usepackage{ mathrsfs }

\newtheorem{theorem}{Theorem}[section]

\newtheorem{definition}[theorem]{Definition}

\newcommand{\TOT}{\Leftrightarrow}

\newcommand{\To}{\Rightarrow}

\date{}

\begin{document}
\thispagestyle{empty}

\title{The Role of the Fifth Postulate in the Euclidean Construction of Parallels}

\author{{\bfseries Iosif Petrakis}\\
   Ludwig-Maximilians Universit\"{a}t M\"{u}nchen, Germany\\
   petrakis@mathematik.uni-muenchen.de\\
} \maketitle

\begin{abstract}
We ascribe to the Euclidean Fifth Postulate a genuine constructive role, which makes it absolutely necessary in the parallel
construction. For that, we present a reconstruction of the general principles underlying the Euclidean construction of a 
geometric property. As a consequence, the epistemological role of Euclidean constructions is revealed. We also examine 
some first implications of our interpretation to the relation between Euclidean and non-Euclidean geometries. 
The Bolyai construction of limiting parallels is also discussed from the Euclidean point of view, as this is reconstructed 
here.
\end{abstract}

\section{The Standard Interpretation of the Fifth Postulate}
\label{sec: SI}

From Proclus up to our days a hermeneutic tradition regarding the Fifth Postulate (FP) has been developed, which we call 
the \textit{Standard Interpretation} (SI). According to it, the Euclidean FP, though originally is formulated differently,
actually 
asserts that through a given point outside a given straight line at most a unique parallel straight line can be drawn to it. 
This formulation, commonly known as Playfair's Axiom (PA), is logically equivalent to the original FP. As the existence of 
a parallel line is independent from PA, the addition of PA establishes the existence of exactly one such parallel.
A clear indication of the predominance of SI is that PA was made the standard form of presenting FP 
in the axiomatic formulations of Euclidean geometry.
In order to describe the shortcomings of SI, we give briefly the Euclidean line of presentation of the parallel
construction in a more formal scheme compatible to our later reconstruction.

\begin{definition}\label{def: par}
If $a, b$ and $c$ are Euclidean coplanar straight lines, we define their following properties: 
$$T(a, b, c) :\TOT  c \ \mbox{falls on} \ a \ \mbox{and} \ b,$$
$$Q_{b}(a) :\TOT  a \ \mbox{is parallel to} \ b,$$ 
$$P_{b, c}(a) :\TOT T(a, b, c) \ \mbox{and} \ c \ \mbox{makes the alternate angles equal to one another}.$$
\end{definition}

\noindent
\textbf{The path to the Euclidean parallel construction} can be described as follows:

\begin{itemize}
 \item Proposition 27 of Book I of the Elements, or \textbf{Proposition I.27}, is a criterion of parallelism, which,
 according to Definition~\ref{def: par}, is the implication $P_{b, c}(a) \To Q_{b}(a)$.
 \item \textbf{Proposition I.28} contains two more criteria of parallelism, reducible to that
 of Proposition I.27.  
 \item \textbf{Proposition I.29}: Let $a, b, c$, such that $T(a, b, c)$, then 
$Q_{b}(a) \To P_{b, c}(a).$\\[1mm]
Hence, in Proposition I.29 the converse to the implication of Proposition I.27 is established. In its proof 
Euclid uses FP\footnote{The original formulation of FP is the following: 
if $T(a, b, c)$ and $c$ makes the interior angles less than two right angles ($2\llcorner$),
then $a, b$, if produced indefinitely, meet on that side on which are the angles less than $2\llcorner$.} for 
the first time.
\item \textbf{Proposition I.30}: If $Q_{b}(a)$ and $Q_{b}(c)$, then $Q_{c}(a)$.\\[1mm]
This proposition, the proof of which uses Proposition I.29, is crucial to SI, as
it proves the uniqueness of the parallel line. This uniqueness result though, is not included in the Elements.
\item \textbf{Proposition I.31}: Construction of a straight line $a$, through a given point $A$ outside line $b$, such 
that $Q_{b}(a)$. The construction consists in the construction of lines $c$ and $a$, such that $P_{b, c}(a)$. 
Then, by Proposition I.27, we also get $Q_{b}(a)$.
\end{itemize}

\noindent
\textbf{The core of SI} can be summarised as follows:
\begin{itemize}
 \item The construction in Proposition I.31 requires only Proposition I.27 and hence it is 
independent from FP. Consequently, it could be placed after Proposition I.27 and before Proposition I.29.

\item Within SI the place of the parallel construction after the first use of FP is explained, although not with absolute
certainty, as an expression of Euclid's need, before giving the construction, to place beyond all doubt the fact that 
only one such parallel can be drawn (see~\cite{He56}, p.~316). Actually, this is an argument of Proclus, as
expressed in his Commentary~\cite{Pr92}, pp. 295--296. If the parallel construction 
was placed right after Proposition I.27, then only the
existence of the parallel line would be established.
\end{itemize}

The independence of FP from
the parallel construction within SI is one of the main reasons why mathematicians, 
before the emergence of non-Euclidean geometries, 
used to consider FP as a theorem rather than as a Postulate.
According to SI, the Euclidean line of presentation certifies 
the existence and the uniqueness of the parallel line, hence the 
``true'' meaning of FP is the uniqueness formulation of the parallel line. 
It is this emphasis of SI on the uniqueness of the parallel which 
pushed it forward as a central characteristic of Euclidean geometry. Gradually, the difference between Euclidean 
geometry and non-Euclidean geometries was identified with the different number of parallels they permit.
 
The uniqueness interpretation though, is, in our view, inadequate. In the first place, there is no explanation 
within SI why Euclid preferred his formulation of FP rather than the uniqueness assumption. Furthermore, the examination 
of the Elements shows that Euclid seems indifferent to questions of uniqueness. In the First Postulate (construction 
of a line segment between two points) there is no reference to the uniqueness of the corresponding 
segment, though it is used in
Proposition I.4 in the form ``two straight lines cannot enclose a space''. The circle under construction in the Third Postulate 
(construction of a circle of any center and radius) is not mentioned to be unique either. The investigation  
of the perpendicular constructions in Propositions I.11 and I.12 reveals the same Euclidean indifference to the
uniqueness of a constructed object.

\section{The basic principles of a Euclidean construction and the constructive role of the Fifth Postulate}
\label{sec: principles}

In this section we present a reconstruction of the general principles underlying the Euclidean construction of a 
geometric property. As a consequence, FP plays a crucial role to the parallel construction of Proposition I.31.
The first three Euclidean postulates have a direct constructive role, as they provide the fundamental elements 
for the subsequent line and circle constructions. In our view, the Fourth and the Fifth Postulate have an indirect, 
though genuine, constructive role. Both are less elementary and they participate in the less elementary 
parallel construction.\\[1mm]
\noindent
\textbf{The constructive role of the Fourth Postulate}: It is used in Proposition I.16 (through Proposition I.15), 
which is necessary in the proof of Proposition I.27. By this line of thought, it participates in the construction 
of Proposition I.31. Also, by the Fourth Postulate, the right angle is a fixed and universal standard, to which
other angles can be compared. In this way FP, treating the $2\llcorner$ as a fixed quantity,
``depends'' on the Fourth Postulate.\\[1mm]
To reveal the constructive character of FP, we need to understand the conceptual requirements of ancient
Greek mathematics regarding the nature of geometric constructions as these are embodied in the Euclidean Elements. 
These requirements are not explicitly found in Euclid, but we propose them as an accurate reconstruction of the 
Euclidean constructive spirit\footnote{A first presentation of our reconstruction of the principles 
underlying the Euclidean constructions is found in~\cite{Pe08}.}.\\[2mm]
\noindent
\textbf{Basic principles of the Euclidean construction $K(a, P)$ of an object $a$
satisfying a geometric property $P$}:\\[2mm]
\noindent
\textbf{K1}: 
In $K(a, P)$ an object $a$ is constructed, satisfying, as accurately as possible, the 
definition of $P$\footnote{The expression ``as accurately as possible'' in K1 will be evident in 
section~\ref{sec: epistem}. K1 can also be found, though not as explicitly as here, in the intuitionistic
literature on the concept of species, or an intuitionistic property. A constructive principle such as K1 can be 
detected in Brouwer's notes. Furthermore, for Griss, 
a species is defined by a property of mathematical objects, but such a property can only have a clear sense after
we have constructed an object which satisfies it (see~\cite{He71}, p.~126). For the construction of species 
see also~\cite{Pe08, Pe10}.
}.\\[1mm]
\noindent
\textbf{K2}: If an object $b$ satisfying a geometric property $R$ is used in construction $K(a, P)$, then the
construction $K(b, R)$ must have already been completed.\\[1mm]
\noindent
\textbf{K3}: If $a$ is a geometric object satisfying $P$ and $Q$ is another geometric property, such that 
whenever $a$ satisfies $P$ it satisfies $Q$, but not the converse i.e.,
$P(a) \To Q(a)$, but not $Q(a) \To P(a)$, then $K(a, Q)$ \textit{cannot} be established through $K(a, P)$.\\[1mm]
\textbf{K4}: If $a$ is a geometric object satisfying $P$ and $Q$ is another geometric property, such that 
whenever $a$ satisfies $P$ it satisfies $Q$, and the converse i.e.,
$P(a) \TOT Q(a)$, 
then $K(a, Q)$ \textit{can} be established through $K(a, P)$, and conversely.

\vspace{2mm}

Principle K2 guarantees that $K(a, P)$ does not have constructive gaps, i.e., all geometric concepts used in
construction $K(a, P)$ are already constructed\footnote{Though K2 is very natural to accept, it is not trivial.
In a sense described in section~\ref{sec: neg}, Bolyai's construction of limiting parallels violates it.}.
Principle K3 is the most crucial to our reconstruction of the role of FP in the Euclidean parallel construction.
It guarantees that the construction of the 
abstract object $a$ satisfying property $Q$ cannot be established through the construction of the less 
general property $P$ i.e., construction $K(a, P)$ respects the generality hierarchy of geometric concepts. 
For example, the construction of an isosceles triangle cannot be established through the construction of an
equilateral triangle, since there are isosceles triangles which are not equilateral\footnote{Euclid uses the
concept of an isosceles triangle in Proposition I.5, without providing first a construction of it, because 
this construction is a simple generalization of the equilateral one (Proposition I.1). Evidently, Euclid 
found no reason to include this, strictly speaking, different, but expected construction.}.
Principle K4 guarantees that whenever properties $P$ and $Q$ are logically equivalent, having the same generality,
they do not differ with respect to construction. K4 is the natural complement to K3 and together
they form  the
core of the Euclidean constructive method.
In order to explain the use of the above set of principles in the parallel construction and their relation
to FP, we define the following notions of construction.

\begin{definition}\label{def: defs}
We call a construction $K(a, P)$ \textit{direct}, if it establishes an object $a$ reproducing exactly 
the definition of $P$. In this case we call $P$ a \textit{finite} property. A geometric property
$Q$ is called \textit{infinite}, if it is impossible to give a direct construction of $Q$\footnote{This impossibility 
is not a logical one, but just a result of the definition of $Q$.}. A construction $K(a, Q)$ is called 
\textit{indirect}, if $K(a, Q)$ establishes an object $a$, which satisfies the definition of $Q$ indirectly,
i.e., through a provably equivalent, finite property $P$. 
\end{definition}

Most of Euclidean constructions are direct. For example, at the end of the perpendicular construction of
Proposition I.12 Euclid restates the definition of the perpendicular line, showing that he has constructed
an object which satisfies completely that very definition. In our terminology, 
the property of perpendicularity is finite. 
On the other hand, the parallel property is infinite. Euclid defined parallel lines 
(Definition 23 of Book 1) as straight lines which, being in the same plane and being produced indefinitely 
in both directions, do not meet one another in either direction. It is impossible to give a direct
construction of a line parallel to a given one, since we cannot reproduce the above definition. The 
infinite character of this definition lies in our mental inability to produce a line indefinitely and 
act as if this ever extended object exists as a whole. Each moment we only know a finite part of the on going line,
from which we cannot infer that every extension of it does not meet the given line. The formation of
the parallel line never ends.\\[1mm]
\noindent
\textbf{Euclidean construction of the infinite parallel property}: Euclid gradually constructed (mainly through 
the Fourth Postulate and Propositions I.16 and I.27) the geometric property $P_{b, c}(a)$, which is a finite 
property. Given a line $b$, we can construct directly lines $c$ and $a$ such that $P_{b, c}(a)$ (actually this
is the construction of Proposition I.31), using only the direct construction of Proposition I.23 (construction
of a rectilinear angle equal to a given one, on a given straight line and at a point on it).\\
The implication $P_{b, c}(a) \To Q_{b}(a)$ is established by Proposition I.27, but it would be a 
violation of K3 if construction $K(a, P_{b, c}(a))$ was considered as construction $K(a, Q_{b}(a))$. 
Construction $K(a, P_{b, c}(a))$ can be considered as construction $K(a, Q_{b}(a))$ only if the converse 
implication $Q_{b}(a) \To P_{b, c}(a)$ is proved. Then, $P$ and $Q$ will have the same generality
and then K4 can be applied. In our view, this is why Euclid ``postponed'' the parallel construction, 
placing it after Proposition I.29, which establishes the converse implication.\\[1mm]
\noindent
\textbf{The constructive role of FP}: The fifth postulate is this (intuitively true) proposition, through which
the implication $Q_{b}(a) \To P_{b, c}(a)$ is shown, and then by K4, construction $K(a, P_{b, c}(a))$
in Proposition I.31 is also construction $K(a, Q_{b}(a))$ of parallels. 

\vspace{2mm}

Euclid used his specific formulation of FP so that the proof of Proposition I.29 requires one only conceptual 
step, reaching his goal in the simplest way. According to our interpretation, Euclid does not postpone the use of
FP as long as possible\footnote{For a recent reference to this long repeated view see~\cite{Ha00}.}, recognizing 
its ``problematic'' nature. On the contrary, he uses it exactly the moment he needs it, revealing in 
this way its function.
In Euclid's Elements, if $P$ is a finite property then $K(a, P)$ is always given through $P$ itself and not through 
an equivalent property $Q$ i.e., K4 is not used in constructions of finite properties. It is used only when
an infinite property $Q$ is to be constructed. Otherwise, its function wouldn't be clear.

The indirect construction of an infinite geometric property is not the only way ancient Greeks used to handle 
an infinite property. If an infinite property $Q$ has no finite equivalent, it may have a special case $F$ 
with a strong finite character accompanying the infinite one. 
Infinite anthyphairesis  i.e., infinite continued fraction $Q$ is an infinite property studied in Book X of the
Elements, which does not have a finite equivalent. Periodic anthyphairesis (periodic continued fraction) $F$ 
is a special case of $Q$, which possesses a strong finite character beside its infinity. Although the sequence 
of the quotients forming the periodic continued fraction never ends (infinity of $F$), its finite period expresses our
knowledge of this sequence (finite character of $F$)\footnote{Ancient Greeks had also found a necessary and 
sufficient condition for an infinite anthyphairesis to be periodic (logos criterion). Its knowledge and its 
importance in Plato's system have been developed in recent times in Negrepontis' program on reconstructing Plato
(see e.g.,~\cite{FN21}).}.
In intuitinisitic terms, periodic anthyphairesis is a strong form of law-likeness.

\section{The epistemological role of Euclidean constructions}
\label{sec: epistem}

Our description of the Euclidean constructive principles also reveals the difference between ``Euclidean
construction and ``Euclidean existence''. Let the existential formula
$$\exists{a}Q(a) :\TOT \mbox{there exists a geometric object} \  a \ \mbox{satisfying geometric property} \ Q.$$
In Euclid's Elements $\exists{a}Q(a)$ is established either by $K(a, Q)$ or by $K(a, P)$, where $P(a) \To Q(a)$
(but not the converse). Euclidean geometry (except Eudoxus' theory of ratios) is a basic paradigm of a 
constructive mathematical theory, since existence of a mathematical object or concept is constructively established.
For example, if the construction of Proposition I.31 was placed right after Proposition I.27, that would only
show the existence of a parallel line. This proof of existence though, does not constitute a construction of the concept 
of parallelism. 


The traditionally accepted independence between FP and the construction of Proposition I.31 is based on the 
identification between $\exists{a}Q(a)$ and $K(a, Q)$\footnote{According to Zeuthen~\cite{Ze96}, the main purpose
of a geometric construction is to provide a proof of existence, so the purpose of the FP is to ensure the 
existence of the intersection point of the non parallel lines. This approach fails to see, in our view, the 
difference between existence and construction.}. For Euclid though, the construction of property $Q$ is generally 
an enterprise larger than the exhibition-construction of a single object satisfying $Q$. The parallel 
construction reflects this fact very clearly. According to our reconstruction, 
$\exists{a}Q(a)$ shows that property $Q$ is not void, that it possesses, in modern terms, an extension. 
On the other hand, $K(a, Q)$ \textit{shows that we have found a way to mentally
grasp property $Q$, fully if $Q$ is finite, as much as possible if $Q$ is infinite}.

Traditionally, the Elements are considered as the original model of the axiomatic method and logical
deduction. In our view, they also constitute a first model of the constructive method, quite different though, 
from modern constructivism in mathematics (see e.g.,~\cite{BIRS23}). It is this 
combination of the axiomatic and the constructive method that reflects the philosophical importance of the
Elements. For the first time in the history of mathematics a mathematical theory answers simultaneously the
ontological and the epistemological problem of the mathematical concepts involved. The ontology of Euclidean
geometric objects and concepts is of mental (and not empirical) nature. Almost certainly Euclidean ontology 
is very close to Platonic ontology\footnote{Euclid was a Platonist and his definitions are closely related to the Platonic 
ones (see~\cite{He71}, p.~168).}.
This mental ontology of mathematical concepts imposes the 
constructive method. It is the construction of mathematical concepts which provides their study with a firm epistemology.
Euclid does not only care about the logical relations between geometric concepts and objects. He also needs 
to answer the main epistemological question: \textit{how do we understand the concepts that we employ in our deductions}? 
And his answer, in our view, is: \textit{we understand them because we construct them}. Thus, \textit{geometric constructions form the 
indispensable epistemology of Euclidean geometry}\footnote{ For a more recent discussion on the role of Euclidean constructions 
see~\cite{Kn83} and~\cite{Ha03}. In our opinion, the interpretations proposed there
are not satisfactory.}.

\section{On the relation between Euclidean and non-Euclidean Geometry}
\label{sec: neg}

It is impossible here to study fully the relation between Euclidean geometry and non-Euclidean geometries. 
We shall only stress some points that can be derived from our previous analysis.
There is a traditional view regarding the above relation too. According to it, the two geometries can be seen as
mathematical structures of the same kind, differing only in the number of parallels they allow. 
One such common mathematical framework is the concept of Hilbert plane\footnote{This framework is not as
absolute as it is often named, since it does
not include the elliptic plane, in which there exist no parallels at all, and every line through the pole of a given 
line is perpendicular to it. Hilbert's classic work~\cite{Hi71} is still the best introduction to Hilbert planes.
A more absolute framework that includes elliptic geometry, is the concept of a Bachmann plane, or metric plane
(see~\cite{Ba71}).}. A Hilbert plane is a system of points, lines and planes satisfying the well known Hilbert
axioms of incidence, betweenness and congruence. In a Hilbert plane the parallel line (as any other geometric property) is not 
constructed, only its existence is established. A general Hilbert plane is neutral with respect to the uniqueness of the parallel line.
A Euclidean plane is a Hilbert plane allowing only one parallel, while a hyperbolic plane is a 
Hilbert plane allowing more than one parallels. 
The consequences of this ``coexistence'' of Euclidean and non-Euclidean geometry were very serious.
The foundations of mathematics and mathematics itself
were influenced immensely from the loss of the a priori character of Euclidean geometry. As 
Euclidean geometry became just one possible geometry, the Kantian a priori suffered a serious blow and especially 
the a priori of space. As a result, all major foundational
programs if mathematics in the twentieth century rest either on a Kantian a priori of discrete nature (e.g., Brouwer's 
intuitionism~\cite{He71} and Bishop's constructivism~\cite{Bi67, BB85}), or on a purely logical substratum 
(e.g., Frege's logicism)\footnote{Putnam's assessment 
in~\cite{Pu75}, p.~x, is characteristic: ``...the overthrow of Euclidean geometry is the most important event in the history of science for
the epistemologist''.}. 

Our reconstruction of the parallel construction suggests a class with this traditional view. In our opinion, 
Euclidean geometry has a constructive character, of a specific type, which non-Euclidean geometry lacks.
Of course, this opinion echoes Kant. In~\cite{We95},
p.~1, Webb remarks the following:
\begin{quote}
[It was a commonplace of older Kantian scholarship that the discovery of non-euclidean geometry undermined his 
theory of the synthetic a priori status of geometry. It is commonplace of newer Kant scholarship that he already knew
about non-euclidean geometry from his friend Lambert, one of the early pioneers of this geometry, and that in fact 
its very possibility only reinforces Kant's doctrine that euclidean geometry is synthetic a priori because only its
concepts are constructible in intuition.]
\end{quote}

The common language of Hilbert planes (or any such common mathematical framework) ignores the role and the necessity of FP 
in the parallel construction, just as the epistemological role of constructions. Modern geometry, developed within classical logic,
seems quite indifferent to epistemological questions. We can only indicate here that the two geometries 
are not directly comparable, from the constructive point of view.
Consequently, Euclidean geometry has not lost its a priori character. 
Next we explain why Bolyai's construction of limiting
parallels is problematic from the Euclidean point of view, as this is reconstructed above.

\vspace{2mm}

\noindent
\textbf{The path to the Bolyai construction of limiting parallels} can be described as follows:

\begin{itemize}
 \item A \textit{hyperbolic plane} is a Hilbert plane satisfying \textbf{Lobachevsky's axiom} (L): If $a$ is a line and $A$ is
 a point outside $a$, there exist rays $Ab, Ac$, not on the same line, which do not intersect $a$, 
 and each ray $Ad$ in the angle $bAc$ intersects $a$.
  
 \item \textbf{Proposition 4.1}: A triangle in a hyperbolic plane has angle sum less than $2\llcorner$.

 \item A quadrilateral $PQRS$ is a \textit{Lambert quadrilateral}, if it has right angles at $P, Q$ and $S$.
 
 \item \textbf{Proposition 4.2}: In a hyperbolic plane the fourth angle (the angle at $R$) of a Lambert quadrilateral $PQRS$
is acute, and a side adjacent to it is greater than its opposite side ($QR > PS$ and $SR > PQ$).

\item \textbf{Proposition 4.3}: Suppose we are given a line $a$ and a point $P$ not on $a$, in a hyperbolic plane. Let $PQ$
be the perpendicular to $a$. Let $m$ be a line through $P$, perpendicular to $PQ$. Choose any point $R$ on $a$, and
let $RS$ be the perpendicular to $m$. If $Pc$ is a limiting parallel ray intersecting $RS$ at $X$, then $PX = QR$.

\item \textbf{Elementary Continuity Principle} (ECP): If one endpoint of a line segment is inside a circle and the other
outside, then the segment intersects the circle.

\item \textbf{Bolyai's construction}: Consider a hyperbolic plane satisfying ECP. 
Suppose we are given a line $a$ and a point $P$ not on $a$. Let $PQ$ be the perpendicular to $a$. Let $m$ be a line
through $P$, perpendicular to $PQ$. Choose any point $R$ on $a$, and let $RS$ be the perpendicular to $m$. Then the 
circle of radius $QR$ around $P$ will meet the segment $RS$ at a point $X$, and the ray $PX$ will be the limiting 
parallel ray to $a$ through $P$.

 \end{itemize}

The proof of Bolyai's construction goes as follows: 
as $Q = \llcorner$ and $PR > QR$, by Proposition 4.1 the angle at $Q$ is the largest angle in the triangle $PQR$. 
Moreover, $PS < QR$, since $PQRS$ is a Lambert quadrilateral satisfying the hypothesis of
Proposition 4.2. Consequently, the endpoints $R$ and 
$S$ of the segment $RS$ are outside and inside the circle $(P, QR)$. By ECP the segment $RS$ intersects $(P, QR)$ at a
(unique) point $X$, and $PX$ is the limiting parallel ray to $a$ through $P$, since L guarantees its existence and 
by Proposition 4.3 it satisfies $PX = QR$. 

The curious feature of this proof is that we prove that this construction works by assuming first 
(via axiom L) that the object we wish to construct already exists. This curiosity is stressed by Hartshorne 
in~\cite{Ha00} p.~398. As this presupposed existence of the limiting parallel is axiomatic and not constructive,
Bolyai's construction violates the Euclidean Principle K2.
Another aspect of the problematic character of Bolyai's construction is related to the constructive principles K3 and K4. 
Proposition 4.3 is in analogy to Proposition I.29, since it can be written in the form 
$$\mbox{L} \To PX = QR.$$
In our language, L is an infinite property and $PX = QR$ is a finite one. In order to consider, from the 
Euclidean point of view, the direct construction of $X$ as the construction of the limiting ray, we have 
to prove directly, in a hyperbolic plane satisfying ECP, the analogue to Proposition I.27:
$$PX = QR \To L.$$
Such a direct proof has not yet been found. Although the above line and circle construction of
the most important concept of hyperbolic geometry shows Bolyai's constructive sensitivity, it does not
satisfy the constructive principles of the Euclidean parallel construction.
The usual proof of the existence of limiting parallel is based on \textit{Dedekind's continuity 
axiom} (D) (see e.g.,~\cite{Gr80} p.~156). According to it, any (set-theoretic) separation of points on
a line i.e., a Dedekind cut, is produced by a unique point. 
Axiom D seems unfathomable from the Euclidean point of view, maybe because of its set-theoretic nature.
The question whether Bolyai's construction could be used to prove the existence of 
a limiting parallel for a system of axioms that includes ECP but not D, was naturally raised by 
Greenberg in~\cite{Gr79a}. Pejas, working in the framework of Bachmann plane geometry, a 
geometry without betweenness and continuity
axioms, succeeded to classify all Hilbert planes\footnote{A Hilbert plane corresponds to an ordered Bachmann 
plane with free mobility. As Greenberg puts it in~\cite{Gr79b}, Hilbert's approach is thus incorporated 
into Klein's Erlangen program, whereby the group of motions becomes the primordial object of interest. 
For Pejas' classification theorem see~\cite{Pe61}.}. 
Greenberg, using Pejas' classification of Hilbert planes,
managed in~\cite{Gr79a} to answer his question positively.\\[1mm]
\textbf{Proposition 4.4} (Pejas-Greenberg): If ECP holds and the fourth angle of a Lambert quadrilateral 
is acute, then Bolyai's construction gives the two lines through $P$ that have a ``common perpendicular 
at infinity'' with $a$ through the ideal points at which they meet $a$. Among Hilbert planes satisfying 
ECP, the Klein models are the only ones which are hyperbolic, and Bolyai's construction gives the asymptotic 
parallels for them.

\vspace{2mm}

An important corollary of Proposition 4.4 is the following proposition.\\[1mm]
\textbf{Proposition 4.5}: Every Archimedean, non-Euclidean\footnote{A Hilbert plane $P$ is called non-Euclidean if 
PA fails in $P$.} 
Hilbert plane in which ECP holds is hyperbolic.

\vspace{2mm}

Though Pejas-Greenberg managed to show that the Bolyai construction does yield the limiting parallel replacing 
D with more elementary axioms, their proof is indirect, since it is based on a classification theorem. 
Hence, from the Euclidean, constructive point of view, there is still no direct constructive proof of the concept 
of limiting parallel. We conjecture that such a proof cannot be found.

%
%
%
%
%
%
%
%
%
%
%

\end{document}